\documentclass[10pt]{article}
\usepackage{amssymb,latexsym,amsmath,epsfig,amsthm}

\usepackage{amsfonts}
\usepackage[dvipsnames]{xcolor}

\makeatletter

\renewcommand\section{\@startsection {section}{1}{\z@}
{-30pt \@plus -1ex \@minus -.2ex}
{2.3ex \@plus.2ex}
{\normalfont\normalsize\bfseries\boldmath}}

\renewcommand\subsection{\@startsection{subsection}{2}{\z@}
{-3.25ex\@plus -1ex \@minus -.2ex}
{1.5ex \@plus .2ex}
{\normalfont\normalsize\bfseries\boldmath}}

\renewcommand{\@seccntformat}[1]{\csname the#1\endcsname. }

\makeatother

\newtheorem{theorem}{Theorem}

\newtheorem{proposition}{Proposition}

\theoremstyle{definition}
\newtheorem{definition}{Definition}
\newtheorem{conjecture}{Conjecture}

\begin{document}

\begin{center}
\uppercase{\bf Will the real Hardy--Ramanujan formula please stand up?}
\vskip 20pt
{\bf Stephen DeSalvo}\\
{\tt stephen.desalvo@gmail.com}
\end{center}
\vskip 20pt
\vskip 30pt

\centerline{\bf Abstract}
\noindent
The Hardy--Ramanujan formula for the number of integer partitions of $n$ is one of the most popular results in partition theory. 
While the unabridged final formula has been celebrated as reflecting the genius of its authors, it has become all too common to attribute either some simplified version of the formula which is not as ingenious, or an alternative more elegant version which was expanded on afterwards by other authors.
We attempt to provide a clear and compelling justification for distinguishing between the various formulas and simplifications, with a summarizing list of key take-aways in the final section. 

\pagestyle{myheadings}
\thispagestyle{empty}
\baselineskip=12.875pt
\vskip 30pt

\section{As simple as possible, but $n$}
Five score and several years ago, G.~H.~Hardy and S.~Ramanujan wrote a consequential paper titled \emph{Asymptotic formul\ae\ in combinatory analysis} \cite{HR}.
In it was an application of a well-known method to the problem of counting \emph{integer partitions}.
An integer partition of size~$n$ is a multiset of \emph{positive} integers whose sum is $n$.  Let $p(n)$ denote the number of integer partitions of size~$n$.  
Hardy and Ramanujan gave their result in a sharp form,~\eqref{realHR} below.  A simple corollary of~\eqref{realHR}, and arguably the most popularly cited formula, called the Hardy--Ramanujan asymptotic for $p(n)$, is given by~\eqref{pn1}. 
We write $f(n) \sim g(n)$ whenever $f(n)/g(n) \to 1$ as $n \to \infty$.
We have 

\begin{equation}
\tag{1.41}
\label{pn1}
p(n) \sim \frac{e^{\sqrt{\frac{2}{3}}\pi \sqrt{n}}}{4\sqrt{3} n} \ .
\end{equation}

While this formula \emph{does} appear in the form above in their original paper, it is too often cited as the sole representation of the ``genius" of its authors.  
This is rather unfortunate since this corollary does not convey any error estimates, and, more importantly, several elementary proofs have since been provided; see for example~\cite{CanfieldElementary} and \cite{ErdosElementary}.  
What's more, the significance of this formula is downplayed in the original paper itself via numerical calculations for $n=10, 20, 50 ,80$, where Hardy and Ramanujan write (with $\varpi$ denoting the right-hand side of \eqref{pn1} and $p$ denoting the left-hand side of \eqref{pn1}), ``It will be observed that the progress of $\varpi/p$ towards its limit unity is not very rapid, and that $\varpi-p$ is always positive and appears to tend rapidly to infinity."

What has not been replicated is an elementary proof of their simplest quantitative version of~\eqref{pn1} above, i.e.,~\eqref{realHR} below, which is only slightly less terse. 
The notation $f(n) = O(g(n))$ means for some finite $C \geq 0$, we have $|f(n)| \leq C\, g(n)$ for all $n$.

\begin{proposition}[Hardy and Ramanujan, 1918~\cite{HR}]\label{propHR}
As $n$ tends to infinity, we have
\begin{equation}\tag{1.55}\label{realHR}
p(n) = \frac{e^{\pi \sqrt{\frac{2}{3}} \sqrt{n - \frac{1}{24}}}}{4\sqrt{3}\left(n-\frac{1}{24}\right)}\left(1 - \frac{1}{\pi \sqrt{\frac{2}{3}} \sqrt{n - \frac{1}{24}}}\right)  + O\left(\frac{e^{\pi \sqrt{\frac{2}{3}}\,  \sqrt{n} / 2}}{n}\right).
\end{equation}
\end{proposition}

Again quoting their original paper, ``The formula (1.55) is an asymptotic formula of a type far more precise than that of (1.41)."
If one were pressed to quote one ``simple" \emph{formula} as distinctly Hardy and Ramanujan's, then the above would suffice, subject to a very strict definition of the meaning of \emph{formula}, which we fully expound upon in Section~\ref{whatisaformula}. 

This may seem like the end of the story, but there is far more to discuss. 
The main theorem of Hardy and Ramanujan, which will be described subsequently, was an asymptotic series with a quite startling etiology and evolution. 
We are principally motivated by a staggering degree of misunderstanding which has seeped its way into both formal expositions on the topic, as well as a plethora of popular culture articles being written which tend to simplify the truth to a degree beyond Einstein's famous quote, ``Everything should be made as simple as possible, but no simpler."
When possible, we will let the relevant authors tell the story in their own words (with references and footnotes locally defined), and attempt to weave one single narrative spanning several decades and published articles. 

\section{The old log theater}

Before we delve into the main results, we start with a result even weaker than~\eqref{pn1}; i.e., an asymptotic for $\log p(n)$.  
Hardy and Ramanujan~\cite{HR} provided a simple asymptotic for $\log p(n)$ as follows.  
\begin{quotation}\color{blue}
$\cdots$ 
the question whether a constant $C$ exists such that
\[ \log p(n) \sim C \sqrt{n}. \]
We prove that this is so in section 3.  Our proof is still, in a sense, ``elementary".   It does not appeal to the theory of analytic functions, depending only on a general arithmetic theorem concerning infinite series; $\cdots$ It shows that 
\begin{equation}\tag{1.35} C = \frac{2\pi}{\sqrt{6}}; \end{equation} 
in other words that
\begin{equation}\tag{1.36}
p(n) = \exp \left\{ \pi \sqrt{\left(\frac{2n}{3}\right)}(1+\epsilon)\right\}, 
\end{equation}
where $\epsilon$ is small when $n$ is large. 
$\cdots$ it is equally possible to prove~(1.36) by reasoning of a more elementary, though more special, character: we have a proof, for example, based on the identity
\[ n p(n) = \sum_{\nu=1}^n \sigma(\nu) p(n-\nu), \]
where $\sigma(\nu)$ is the sum of the divisors of $\nu$, and a process of induction.
But we are at present unable to obtain, by any method which does not depend upon Cauchy's theorem, a result as precise as that which we state in the next paragraph
$\cdots$.
\end{quotation}

The result referred to in the ``next paragraph" is Equation~\eqref{pn1}.  Skipping ahead to Section 3 of~\cite{HR}:  

\begin{quotation}\color{blue}
3.1.  The value of the constant
\[ C = \lim \frac{\log p(n)}{\sqrt{n}}, \]
is most naturally determined by the use of the following theorem.

\emph{If $g(x) = \sum a_n x^n$ is a power-series with positive coefficients, and }
\[ \log g(x) \sim \frac{A}{1-x}\]
\emph{when $x \to 1$, then} 
\[ \log s_n = \log (a_0 + a_1 + \ldots + a_n) \sim 2\sqrt{A\, n} \]
\emph{when $n \to \infty$.}

This theorem is a special case\footnote{\color{blue} L.c.~p.~129 (with $\alpha=1$)} of Theorem~C in our paper already referred to.

Now suppose that
\[ g(x) = (1-x) f(x) = \sum \{p(n) - p(n-1)\} x^n = \frac{1}{(1-x^2)(1-x^3)(1-x^4)\ldots} \ .\]
Then
\[ a_n = p(n) - p(n-1)\]
is plainly positive.  And
\begin{equation}\tag{3.11} \log g(x) = \sum_{2}^\infty \log \frac{1}{1-x^\mu} = \sum_1^\infty \frac{1}{\nu} \frac{x^{2\nu}}{1-x^\nu} \sim \frac{1}{1-x} \sum_1^\infty \frac{1}{\nu^2} = \frac{\pi^2}{6(1-x)},\end{equation}
when $x\to 1.$\footnote{\color{blue} This is a special case of much more general theorems : see $\cdots$}
Hence
\begin{equation}\tag{3.12} \log p(n) = a_0 + a_1 + \ldots + a_n \sim C \sqrt{n},\end{equation}
where $C = 2\pi /\sqrt{6} = \pi \sqrt{\frac{2}{3}}$ \ldots \ .
\end{quotation}

P.~Erd\H{o}s~\cite{ErdosElementary} later explicitly wrote out the elementary proof referred to above, as the first step of an elementary proof of~\eqref{pn1}: 

\begin{quotation}\color{RoyalPurple}
\[\cdots \qquad c = \pi \left(\frac{2}{3}\right)^{\frac{1}{2}} \qquad \cdots \]
The starting point will be the following identity:
\begin{equation}\tag{E2} n p(n) = \sum_{v=1} \sum_{k=1} v p(n-kv), \quad p(0) = p(-m) = 0. \end{equation}

(We easily obtain (E2) by adding up all the $p(n)$ partitions of $n$, and noting that $v$ occurs in $p(n-v)$ partitions.). (E2) is of course well known.  In fact, Hardy and Ramanujan state in their paper\footnote{\color{RoyalPurple} Hardy, Ramanujan, Asymptotic formulae in combinatory analysis, Proc. London Math. Soc. 17, (1918), pp.~75-115} that by using (E2) they have obtained an elementary proof of 
\begin{equation}\tag{E3} \log p(n) \sim c n^{\frac{1}{2}}. \end{equation}
The proof of (E3) is indeed easy.  First we show that 
\begin{equation}\tag{E4} p(n) < e^{c\, n^{\frac{1}{2}}}. \end{equation}

We use induction.  (E4) clearly holds for $n=1$.  By (E2) and the induction hypothesis we have
\[ n p(n) < \underset{kv < n}{\sum_{v=1} \sum_{k=1}} ve^{c(n-kv)^{\frac{1}{2}}} < \sum_{v=1}^\infty \sum_{k=1}^{\infty} v e^{c n^{\frac{1}{2}} - ckv/2n^{\frac{1}{2}}} = e^{cn^{\frac{1}{2}}} \sum_{k=1}^\infty  \frac{e^{-kc/2n^{\frac{1}{2}}}}{(1-e^{kc/2n^\frac{1}{2}})^2}\, .\]
Now it is easy to see that for all real $x$, $\frac{e^{-x}}{(1-e^{-x})^2} < \frac{1}{x^2}$.  Thus
\[ n p(n) < e^{c\, n^{\frac{1}{2}}} \sum_{k=1}^\infty \frac{4n}{c^2 k^2} = n e^{c\, n^{\frac{1}{2}}}, \]
which proves~(E4). 

\noindent Similarly but with slightly longer calculations, we can prove that for every $\epsilon>0$ there exists an $A > 0$ such that
\begin{equation}\tag{E5} p(n) > \frac{1}{A} e^{(c-\epsilon)n^{\frac{1}{2}}}. \end{equation}
(E4) and (E5) clearly imply (E3).
\end{quotation}
 
We say $f(n) = o(g(n))$ for positive function $g$ whenever $f(n) / g(n) \to 0$ as $n\to\infty$. 
In particular, $f(n) = o(1)$ implies that $f(n)$ tends to zero as $n \to \infty$.  
To summarize, the above arguments establish 
\[ \log p(n) = \sqrt{\frac{2}{3}}\, \pi \sqrt{n} + o(\sqrt{n}). \]
A restatement of~\eqref{pn1} is 
\begin{equation*}
\log p(n) = \sqrt{\frac{2}{3}}\, \pi\, \sqrt{n} - \log(n) - \log(4\sqrt{3}) + o(1).
\end{equation*}
Erd\H{o}s established in~\cite{ErdosElementary}, through completely elementary means, that there exists an $a > 0$ such that
\begin{equation*}
\log p(n) = \sqrt{\frac{2}{3}}\, \pi\, \sqrt{n} - \log(n) - \log(a) + o(1).
\end{equation*}
A restatement of~\eqref{realHR} is 
\begin{align*}\label{logHR}
\log p(n) & = \sqrt{\frac{2}{3}}\, \pi\, \sqrt{n-\frac{1}{24}} - \log\left( n-\frac{1}{24}\right) - \log(4\sqrt{3}) \\
 &  + \log\left(1 - \frac{1}{\pi \sqrt{\frac{2}{3}} \sqrt{n - \frac{1}{24}}}\right)  + O\left(e^{-\pi \sqrt{\frac{2}{3}}\, \sqrt{n}/2}\right).
\end{align*}

\section{Quite startling}

The generating function of $p(n)$ is easily seen to be
\begin{equation}\label{generating function}
f(x) = \prod_{i=1}^\infty (1-x^i)^{-1}.
\end{equation}
Quoting~\cite{HR} on the origins of their approach:

\begin{quotation}\color{blue}
The function $p(n)$ may, of course, be expressed in the form of an integral
\begin{equation}\tag{1.21} p(n) = \frac{1}{2\pi i} \int_{\Gamma} \frac{f(x)}{x^{n+1}}\, dx, \end{equation}
by means of Cauchy's theorem, the path $\Gamma$ enclosing the origin and lying entirely inside the unit circle.  The idea which dominates this paper is that of obtaining asymptotic formul\ae\ for $p(n)$ by a detailed study of the integral~(1.21).  This idea is an extremely obvious one; it is the idea which has dominated nine-tenths of modern research in the analytic theory of numbers: and it may seem very strange that it should never have been applied to this particular problem before. Of this there are no doubt two explanations. The first is that the theory of partitions has received its most important developments, since its foundation by Euler, at the hands of a series of mathematicians whose interests have lain primarily in algebra. The second and more fundamental reason is to be found in the extreme complexity of the behaviour of the generating function $f(x)$ near a point of the unit circle.

It is instructive to contrast this problem with the corresponding problems which arise for the arithmetical functions $\pi(n)$, $\zeta(n)$, $\psi(n)$, $\mu(n)$, $d(n)$, $\ldots$ which have their genesis in Riemann's Zeta-function and the functions allied to it.  $\cdots$ to write down the dominant terms involves, as a rule, no difficulty more formidable than that of deforming a path of integration over a pole of the subject of integration and calculating the corresponding residue. 

In the theory of partitions, on the other hand, we are dealing with functions which do not exist at all outside the unit circle. Every point of the circle is an essential singularity of the function, and no part of the contour of integration can be deformed in such a manner as to make its contribution obviously negligible. Every element of the contour requires special study; and there is no obvious method of writing down a ``dominant term."

 The difficulties of the problem appear then, at first sight, to be very serious. We possess, however, in the formul\ae\ of the theory of the linear transformation of the elliptic functions, an extremely powerful analytical weapon by means of which we can study the behavior of $f(x)$ near any assigned point of the unit circle\footnote{\color{blue}See G.~H.~Hardy and J.~E.~Littlewood, ``Some problems of Diophantine approximation (II: The trigonometrical series associated with the elliptic Theta-functions)", \emph{Acta Mathematica,} Vol.~37, 1914, pp.~193--238, for applications of the formul\ae\ to different but not unrelated problems.}. It is to an appropriate use of these formul\ae\ that the accuracy of our final results, an accuracy which will, we think, be found to be quite startling, is due.
\end{quotation}

A modern use of the words ``quite startling" in a mathematics paper would be considered 
inappropriate and exceedingly self-important.  On the other hand, when considering both the context of its early 1900s origins as well as the depth and accuracy of their final result, we hope to reaffirm to the reader of the appropriateness of that remark.
 
\section{How low can you go?}

Again quoting from Hardy and Ramanujan's original paper~\cite{HR}:

\begin{quotation}\color{blue}
\begin{equation}\tag{1.53}
C = 2\pi / \sqrt{6} = \pi \sqrt{\frac{2}{3}}, \qquad \lambda_n = \sqrt{n-\frac{1}{24}}.
\end{equation}

 ... \\

The next step is naturally to direct our attention to the singular point of $f(x)$ next in
importance after that at $x = 1$, viz., that at $x = -1$ $\cdots$ 
No new difficulty of principle is involved, and we find that
\begin{equation}\tag{1.61}
p(n) = \frac{1}{2\pi \sqrt{2}} \frac{d}{dn} \left(\frac{e^{C \lambda_n}}{\lambda_n}\right) + \frac{(-1)^n}{2\pi} \frac{d}{dn} \left(\frac{e^{\frac{1}{2}C \lambda_n}}{\lambda_n}\right) + O(e^{D\sqrt{n}}),
\end{equation}
where $D$ is now any number greater than $\frac{1}{3} C$. \\
$\cdots$ \\
The next two terms in the approximate formula are found to be
\begin{equation*}
\frac{\sqrt{3}}{\pi\sqrt{2}} \cos\left(\frac{2}{3}\,n\,\pi - \frac{1}{18}\pi\right) \frac{d}{dn} \left(\frac{e^{\frac{1}{3}C\lambda_n}}{\lambda_n}\right)
\end{equation*}
and
\begin{equation*}
\frac{\sqrt{2}}{\pi} \cos\left(\frac{1}{2}\,n\,\pi - \frac{1}{8}\pi\right) \frac{d}{dn} \left(\frac{e^{\frac{1}{4}C\lambda_n}}{\lambda_n}\right).
\end{equation*}

As we proceed further, the complexity of the calculations increases.  \\
$\cdots$ \\
But it is plain that, by taking a sufficient number of terms, we can find a formula in which the error is
\begin{equation*}
O(e^{C\lambda_n/\nu}),
\end{equation*}
where $\nu$ is a fixed but arbitrarily large integer.
\end{quotation}

The above formulas demonstrate how they are able to apply their method iteratively to produce an asymptotic series expansion for $p(n)$, with ever decreasing Big-O error terms.  As can be seen, each of these terms is itself simple and explicit. 

\begin{quotation}\color{blue}1.7.\  A final question remains. We have still the recourse of making $\nu$ a function of $n$, that
is to say of making the number of terms in our approximate formula itself a function of $n$.
In this way we may reasonably hope, at any rate, to find a formula in which the error is
of order less than that of any exponential of the type $e^{an}$; of the order of a power of $n$, for example, or even bounded.
\end{quotation}

This is perhaps the single greatest setup for what turns out to be a quite startlingly accurate final formula.  Not only are they able to take $\nu$ as a function of $n$, they are able to prove a quantitative error rate that tends to zero!  

\begin{theorem}[Hardy and Ramanujan, 1918~\cite{HR}]
\label{HRfull}
Suppose that
\begin{equation}\tag{1.71}
\phi_q(n) = \frac{\sqrt{q}}{2\pi\sqrt{2}} \frac{d}{dn}\left(\frac{e^{C \lambda_n/q}}{\lambda_n}\right),
\end{equation}
where $C$ and $\lambda_n$ are defined by the equations~$(1.53)$, for all positive integral values of $q$; that $p$ is a positive integer less than and prime to $q$; that $\omega_{p,q}$ is a $24q$-$th$ root of unity, defined when $p$ is odd by the formula
\begin{equation}\tag{1.721}
\omega_{p,q} = \left(\frac{-q}{p}\right)\exp\left[-\left\{\frac{1}{4}(2-pq-p)+\frac{1}{12}\left(q-\frac{1}{q}\right)\left(2p-p'+p^2p'\right)\right\}\pi i\right],
\end{equation}
and when $q$ is odd by the formula
\begin{equation}\tag{1.722}
\omega_{p,q} = \left(\frac{-p}{q}\right)\exp\left[-\left\{\frac{1}{4}(q-1)+\frac{1}{12}\left(q-\frac{1}{q}\right)\left(2p-p'+p^2p'\right)\right\}\pi i\right],
\end{equation}
where $(a/b)$ is the symbol of Legendre and Jacobi\footnote{Quoting the footnote from~\cite{HR}: {\color{blue} See Tannery and Molk, i.e., pp. 104--106, for a complete set of rules for the calculation of the value of $(a/b)$, which is, of course, always $1$ or $-1$.  When both $p$ and $q$ are odd it is indifferent which formula is adopted.}}, and $p'$ is any positive integer such that $1+pp'$ is divisible by $q$; that
\begin{equation}\tag{1.73}\label{Aqn}
A_q(n) = \sum_{(p)} \omega_{p,q} e^{-2np\pi i/q};
\end{equation}
and that $a$ is any positive constant, and $\nu$ the integral part of $a\sqrt{n}$.

Then
\begin{equation}\tag{1.74}\label{HRsum}
p(n) = \sum_1^\nu A_q \phi_q + O(n^{-1/4}),
\end{equation}
so that $p(n)$ is, for all sufficiently large values of $n$, the integer nearest to
\begin{equation}\tag{1.75}\label{HRFormula}
\sum_{1}^\nu A_q \phi_q.
\end{equation}
\end{theorem}

It is easy to understand why their main theorem is rarely quoted, as it requires much more than a passing glance to effectively parse.  On the other hand, Equation~\eqref{realHR} is arguably just as concrete and simple as Equation~\eqref{pn1}, and offers a glimpse into the potential depths of their results without being obtuse.

The insight leading to this remarkable result is articulated by J.~E.~Littlewood; we refer to the following quote from the book \emph{The Theory of Partitions} \cite{Andrews} by George Andrews:

\begin{quotation}\color{Emerald}
The story of the Hardy and Ramanujan collaboration on this formula is an amazing one, and is perhaps best told by J.~E.~Littlewood in his fascinating review of the \emph{Collected Papers of Srinivasa Ramanujan} in the \emph{Mathematical Gazette}, Vol.~14 (1929, pp.~427--428): \\

I must say something finally of the paper on partitions \ldots 

The reader does not need to be told that this is a very astonishing theorem, and he will readily believe that the methods by which it was established involve a new and important principle, which has been found very fruitful in other fields.  The story of the theorem is a romantic one.  (To do it justice I must infringe a little the rules about collaboration.  I therefore add that Prof. Hardy confirms and permits my statements of bare fact.)  One of Ramanujan's Indian conjectures was that the first term of~\eqref{HRFormula} was a very good approximation to $p(n)$; this was established without great difficulty.  At this stage the $n-(1/24)$ was represented by a plain $n$ -- the distinction is irrelevant.  From this point the real attack begins.  The next step in the development, not a very great one, was to treat~\eqref{HRFormula} as an ``asymptotic" series, of which a fixed number of terms (e.g.~$\nu=4$) were to be taken, the error being of the order of the next term.  But from now to the very end Ramanujan always insisted that much more was true than had been established: ``there must be a formula with error $O(1)$."  This was his most important contribution; it was both absolutely essential and most extraordinary.  A severe numerical test was now made, which elicited the astonishing facts about $p(100)$ and $p(200)$.  Then $\nu$ was made a function of $n$; this \emph{was} a very great step, and involved new and deep function-theory methods that Ramanujan obviously could not have discovered by himself.  The complete theorem thus emerged.  But the solution of the final difficulty was probably impossible without one more contribution from Ramanujan, this time a perfectly characteristic one.  As if its analytical difficulties were not enough, the theorem was entrenched also behind almost impregnable defences of a purely formal kind.  The form of the function~$\phi_q(n)$ is a kind of indivisible unit; among many asymptotically equivalent forms it is essential to select exactly the right one.  Unless this is done at the outset, and the $-1/24$ (to say nothing of the $d/dn$) is an extraordinary stroke of formal genius, the complete result can never come into the picture at all.  There is, indeed, a touch of real mystery.  If only we \emph{knew} there was a formula with error $O(1)$, we might be forced, by slow stages, to the correct form of $\phi_q$.  But why was Ramanujan so certain there \emph{was} one?  \emph{Theoretical} insight, to be the explanation, had to be of an order hardly to be credited.  Yet it is hard to see what numerical instances could have been available to suggest so strong a result.  And unless the form of $\phi_q$ was known already, \emph{no} numerical evidence could suggest anything of the kind -- there seems no escape, at least, from the conclusion that the discovery of the correct form was a single stroke of insight.  We owe the theorem to a singularly happy collaboration of two men, of quite unlike gifts, in which each contributed the best, most characteristic, and most fortunate work that was in him.  Ramanujan's genius did have this one opportunity worthy of it.
\end{quotation}

\section{To infinity, and beyond!}
\subsection{Of the greatest interest}

Despite the surprisingly accurate final results, there are some very important caveats that Hardy and Ramanujan were themselves very much aware of.  We next quote the additional remarks section of~\cite{HR}:

\begin{quotation}\color{blue}
6.1. The theorem which we have proved gives information about $p(n)$ which is in some
ways extraordinarily exact. We are for this reason the more anxious to point out explicitly
two respects in which the results of our analysis are incomplete.

6.21. We have proved that
\[p(n) = \sum A_q \phi_q + O(n^{-\frac{1}{4}}),\]
where the summation extends over the values of $q$ specified in the theorem, for every fixed
value of $\alpha$; that is to say that, when $\alpha$ is given, a number $K = K(\alpha)$ can be found such
that
\[|p(n) - \sum A_q \phi_q| < K n^{-\frac{1}{4}}\]
for every value of $n$. It follows that
\begin{equation}\tag{6.211}\label{6.211}
p(n) = \left\{\sum A_q \phi_q\right\},
\end{equation}
where $\{x\}$ denotes the integer nearest to $x$, for $n \geq n_0$, where $n_0 = n_0(\alpha)$ is a certain
function of $\alpha$.
\end{quotation}

From a computational point of view this theorem is the greatest possible tease.  
It provides the existence of a finite series of roughly $\sqrt{n}$ terms which gives an exact value of $p(n)$ as long as $n$ is large enough, without being able to specify what actual value of $n$ is large enough!

Another popular misunderstanding is that Hardy and Ramanujan proved a convergent asymptotic series for the exact value of $p(n)$.  
We again quote the additional remarks section of~\cite{HR}:

\begin{quotation}\color{blue}
6.22. The second point of incompleteness of our results is of much greater interest and
importance. We have not proved either that the series
\begin{equation*}
\sum_{1}^\infty A_q \phi_q
\end{equation*}
is convergent, or that, if it is convergent, it represents $p(n)$. Nor does it seem likely that
our method is one intrinsically capable of proving these results, if they are true---a point
on which we are not prepared to express any definite opinion.
\end{quotation}

This remark is a warning to all future readers to be wary of inference beyond proof. 
Despite the amazing accuracy of the formula, along with its ability to get within $\frac{1}{2}$ of the actual value of $p(n)$, there is no guarantee that the series itself converges.
Once again quoting Hardy and Ramanujan~\cite{HR}: 
\begin{quotation}\color{blue}
The question remains whether we can, by an appropriate choice of $\alpha$, secure the truth of \eqref{6.211} for \emph{all} values of $n$, and not merely for all sufficiently large values.  Our opinion is that this is possible, and that it could be proved to be possible without any fundamental change in our analysis.  Such a proof would however involve a very careful revision of our argument.  It would be necessary to replace all formul\ae\ involving $O$'s by inequalities, containing only numbers expressed explicitly as functions of the various parameters employed.  
This process would certainly add very considerably to the length and the complexity of our argument. It
is, as it stands, sufficient to prove what is, from our point of view, of the greatest interest;
and we have not thought it worth while to elaborate it further.
\end{quotation}

Fortunately, for those interested in computational methods and convergence guarantees, this is not the end of the story.
This is where one has to look beyond the original paper to get more definitive answers and context.

\subsection{Useful for numerical computation}
Another consequential part of the story was provided by H.~Rademacher in 1937, announced and summarized in \emph{A convergent series for the partition function $p(n)$} \cite{rademacher1937convergent}, with details in \emph{On the partition function $p(n)$} \cite{Rademacher}.  

\begin{theorem}[Rademacher, 1937~\cite{rademacher1937convergent}]
We have
\begin{equation}\tag{R4}
p(n) = \frac{1}{\pi \sqrt{2}}\ \sum_{k=1}^\infty  A_k(n) \sqrt{k}\, \frac{d}{dn}\left(\frac{\sinh\left(\frac{\pi \sqrt{\frac{2}{3}} \sqrt{n-\frac{1}{24}}}{k}\right)}{\sqrt{n-\frac{1}{24}}}\right),
\end{equation}
where $A_k(n)$ is given by Equation~\eqref{Aqn}.
\end{theorem}

While very similar to Hardy and Ramanujan's formula, the summands are distinct; one replaces $e^{C \lambda_n/q}$ with $2 \sinh(C \lambda_n/q)$.  The most distinguishing aspect of this work is that Rademacher was able to explicitly bound the error of the approximation after any finite number of terms, both for his series and Hardy and Ramanujan's.  It is therefore a computationalist's dream result, as it allows one to round the partial sum after the approximation error has been \emph{proven} to be below $\frac{1}{2}$.
Quoting from~\cite{rademacher1937convergent}:

\begin{quotation}\color{PineGreen}
From (R4) it is possible to derive a definite expression for the $O$-term in (1.74).  We obtain

\begin{align}
\tag{R12} p(n) = &\  \frac{1}{2\pi \sqrt{2}} \sum_{k=1}^N A_k(n) \sqrt{k} \frac{d}{dn}\left(\frac{e^{C\sqrt{n - \frac{1}{24}}/k}}{\sqrt{n-\frac{1}{24}}} \right) \\
\nonumber &+ \Theta \cdot \left\{ \frac{44\pi^2}{225\sqrt{3}}\, N^{-\frac{1}{2}} + \frac{\pi\sqrt{2}}{75} \left(\frac{N}{n-1}\right)^{1/2} \sinh\frac{C\sqrt{n}}{N} \right. \\
\nonumber &\left. +\frac{1}{2\sqrt{3}} \frac{(N+1)^{3/2}}{n-1} e^{-\frac{C\sqrt{n-1}}{N}}\left(\frac{1}{3} +\frac{\sqrt{3}}{5\pi\sqrt{2}} \frac{N+1}{(n-1)^{1/2}}\right)\right\}, \qquad |\Theta| < 1,
\end{align}
a formula which proves to be useful for numerical computation of $p(n)$.
\end{quotation}

Thus, Rademacher opens up the use of Hardy and Ramanujan's formula for precise numerical computation, and provides his own convergent series for $p(n)$ with similar numerical guarantees.  
Hardy and Ramanujan certainly pioneered not just the approach, but also the idea that it could provide exact results after a finite number of terms.
Rademacher, building off of their approach, found a related (though notably distinct) convergent series for $p(n)$ which allowed him to formalize precisely how many of their terms were sufficient.
Thus, while the two series are inexorably intertwined, we believe the notion of a ``Hardy-Ramanujan-Rademacher" formula or series should be avoided altogether.

Note also that at this point it is still not apparent whether Hardy and Ramanujan's infinite series converges.
Fortunately, there is still another impactful player which, in some sense, completes our story.

\subsection{Spinning infinity}
Shortly after Rademacher, D.~H.~Lehmer investigated the numerical properties of both series in a collection of several papers~\cite{Lehmer1, Lehmer2, Lehmer}.  One punchline is the following:

\begin{theorem}[Lehmer, 1937~\cite{Lehmer1}]
Hardy and Ramanujan's infinite series is divergent.
\end{theorem}

With this result one is forced to take a side on what constitutes ``of the greatest interest."  Hardy and Ramanujan's original paper was both novel in its results as well as substantially groundbreaking, even though it was devoid of any concrete numerical guarantees for particular values of $n$.  Rademacher's convergent series, along with the explicit bounds on the partial sums of both series, provides the numerical guarantees for exact computing.
Hardy and Ramanujan's formula also consequently provides an excellent example of the usefulness of divergent series. 

Lehmer proved a great deal regarding the accuracy of the formulas via the terms $A_k(n)$.  He showed that the error in approximating $p(n)$ by $\alpha \sqrt{n}$ terms of either series is $O(n^{-1/2} \log(n))$.
A particularly useful corollary of this work is the following. 

\begin{theorem}[Lehmer, 1939~\cite{Lehmer2}]\label{Lehmer23Theorem}
If only $2n^{1/2}/3$ terms of the Hardy-Ramanujan series (1.74) be taken, the resulting sum will differ from $p(n)$ by less than $1/2$, provided $n > 600.$
\end{theorem}

Let us reflect on this for a moment.  
Using the first $\frac{2}{3} \sqrt{n}$ terms of a \emph{divergent} series, we approach $p(n)$ to within rounding error, and if we continue adding terms ad infinitum we diverge from $p(n)$.  What's more, Lehmer writes in~\cite{Lehmer2}:

\begin{quotation}\color{Mahogany}
The factor of $2/3$ of Theorem~\ref{Lehmer23Theorem} may be made smaller by allowing the lower limit of $n$ to increase.  For example if we wish to take only $n^{1/2}$/2 terms of the series we may do so provided $n > 3600$.  By making a general argument we may easily prove the following:

\begin{theorem}
Let $\delta > 1$ and let $c = \pi (2/3)^{1/2} = 2.565 \cdots$.  Then  $p(n)$ is the nearest integer to the sum of the first $n^{1/2}/\delta$ terms of the Hardy-Ramanujan series provided 
\begin{equation*}
n > \frac{27^{1/2} c^6}{\delta^2} \left\{\frac{\sinh(c\delta)}{c^3\delta^3} + \frac{1}{6}\right\}^3 = O(e^{3c\delta} \delta^{-11}).
\end{equation*}
\end{theorem}

\end{quotation}

The introduction of his paper \emph{On the remainders and convergence of the series for the partition function}~\cite{Lehmer}, beautifully summarizes the main conclusions of his work with respect to integer partitions:

\begin{quotation}\color{Mahogany}
1.  Introduction.  The two series under discussion are
\begin{align}
\tag{L1} p(n) = & \frac{12^{1/2}}{24n-1} \sum_{k=1}^N A_k^\ast(n)\left(1-\frac{k}{\mu}\right) e^{\mu/k} + R_1(n, N), \\
\tag{L2} p(n) = & \frac{12^{1/2}}{24n-1} \sum_{k=1}^N A_k^\ast(n)\left\{\left(1-\frac{k}{\mu}\right) e^{\mu/k} + \left(1+\frac{k}{\mu}\right) e^{\mu/k}\right\} + R_2(n, N),
\end{align}
due respectively to Hardy and Ramanujan~\cite{HR} (1917) and to Rademacher~\cite{Rademacher} (1937).  Here we have introduced the abbreviation
\begin{equation}\tag{L3}
\mu = \mu(n) = (\pi/6)(24n-1)^{1/2} = O(n^{1/2}).
\end{equation}
The coefficients $A^\ast$ are real numbers defined by
\begin{equation}\tag{L4}
A_k^\ast(n) = k^{-1/2}A_k(n),
\end{equation}
where $A_k(n)$ is a complicated sum of $24k$th roots of unity.  The remainders $R_1(n, N)$ and $R_2(n, N)$ are defined by~(L1) and (L2) in which $p(n)$ denotes the number of unrestricted partitions of $n$.

The fact of primary importance about (L2) is that
\begin{equation}\tag{L5}
\lim_{N\to\infty} R_2(n, N) = 0;
\end{equation}
that is to say, the series in (L2) as $N\to\infty$ converges for all $n$ to $p(n)$.  Concerning $R_1(n, N)$ Hardy and Ramanujan proved that for every $\alpha>0$
\begin{equation}\tag{L6}
R_1(n, \alpha n^{1/2}) = O(n^{-1/4}).
\end{equation}
Rademacher~\cite{Rademacher} gave the following estimate for $R_2(n, N)$ in general:
\begin{equation}\tag{L7}
|R_2(n, N)| < \frac{44\pi^2}{225\cdot 3^{1/2}} N^{-1/2} + \frac{\pi\cdot 2^{1/2}}{75} \frac{N^{1/2}}{(n-1)^{1/2}} \sinh \frac{\pi(2n/3)^{1/2}}{N}
\end{equation}
and a more complicated estimate for $R_1(n, N)$ from which (L6) follows in case $N = \alpha n^{1/2}$.  These estimates for the possible errors in (L1) and (L2) permitted for the first time the use of either (L1) or (L2) with absolute assurance.  Using the estimate
\begin{equation}\tag{L8}
|A_k(n) | < 2k^{5/6}
\end{equation}
instead of the trivial
\begin{equation}\tag{L9}
|A_k(n)| < k
\end{equation}
previously employed, the writer obtained~\cite{Lehmer1, Lehmer2}
\begin{equation}\tag{L10}
|R_2(n, N)| < \frac{\pi^2}{3^{1/2}} N^{-2/3} \left\{\left(\frac{N}{\mu}\right)^3 \sinh \frac{\mu}{N} + \frac{1}{6} - \left(\frac{N}{\mu}\right)^2 \right\},
\end{equation}
\begin{align}\tag{L11}
|R_1(n, N)| < \frac{\pi^2N^{7/3}}{3^{1/2}\mu^3} \left\{ \sinh \frac{\mu}{N}  + \right.& \frac{1}{6}\left(\frac{\mu}{N}\right)^3  \\
\nonumber & \left.+ \left(1 + \frac{N}{\mu}\right)\left(\frac{1}{7} + \frac{1}{3} \mu^{1/3}N^{-5/3}\right) \right\}.
\end{align}
If in (L10) and (L11) we substitute $N = \alpha n^{1/2}$, we find that in either case
\begin{equation}\tag{L12}
R_i(n, \alpha n^{1/2}) = O(n^{-1/3})\qquad \qquad (i=1,2).
\end{equation}
In \S2 we show by a simple asymptotic argument that 
\begin{equation}\tag{L13}
R_i(n, \alpha n^{1/2}) = O(n^{-1/2} \log n)\qquad \qquad (i=1,2),
\end{equation}
a result, which in a sense, is best possible.

$\cdots$ \\

Hardy and Ramanujan \cite[p. 107]{HR} 
raised the question of the boundedness of $A_k(n)$ in discussing the possible convergence of (L1) as $N\to\infty$.  In proving the divergence of (L1) the writer \cite{Lehmer1} employed a sequence of $A$'s which, if they tended to zero, did not do so rapidly enough to render (L1) convergent.  Although this showed, in other words, that $R_1(n, N)$ tends to zero for no value of $n$, it did not remove the possibility of $R_1(n, N)$ ultimately oscillating between fixed limits.  

$\cdots$ \\

\noindent From this result it follows that $R_1(n, N)$ does not oscillate between fixed limits, the terms of the series in (L1) being unbounded. 
It follows also that the $k$th term of (L2) is greater in absolute value than $1/k^2$ for an infinity of $k$'s despite the apparent rapidity of its convergence.

\end{quotation}

\subsection{Useful for numerical computation, reprise}

Venturing into the final pages of~\cite{HR}, one sees in Table II the following representation of the terms $A_k(n)$ for $k = 1, 2, \ldots, 18$:
\begin{quotation}\color{blue}
\noindent $A_1 = 1.$ \\
$A_2 = \cos n\pi.$ \\
$A_3 = 2\, \cos(\frac{2}{3} n\, \pi - \frac{1}{18}\pi).$ \\
$A_4 = 2\, \cos(\frac{1}{2} n\, \pi - \frac{1}{8}\pi).$ \\
$A_5 = 2\, \cos(\frac{2}{5} n\, \pi - \frac{1}{5}\pi) + 2 \cos \frac{4}{5}n\pi.$ \\
$A_6 = 2\, \cos(\frac{1}{3} n\, \pi - \frac{5}{18}\pi)$. \\
$A_7 = 2\, \cos(\frac{2}{7} n\, \pi - \frac{5}{14}\pi) + 2\, \cos(\frac{4}{7} n\, \pi - \frac{1}{14}\pi) + 2\, \cos(\frac{6}{7} n\, \pi + \frac{1}{14}\pi)$. \\
$A_8 = 2\, \cos(\frac{1}{4} n\, \pi - \frac{7}{16}\pi) + 2\, \cos(\frac{3}{4} n\, \pi - \frac{1}{16}\pi)$. \\
$A_9 = 2\, \cos(\frac{2}{9} n\, \pi - \frac{14}{27}\pi) + 2\, \cos(\frac{4}{9} n\, \pi - \frac{4}{27}\pi) + 2\, \cos(\frac{8}{9} n\, \pi + \frac{4}{27}\pi)$. \\
$A_{10} = 2\, \cos(\frac{1}{5} n\, \pi - \frac{3}{5}\pi) + 2\, \cos(\frac{3}{5} n\, \pi)$. \\
$A_{11} = 2\, \cos(\frac{2}{11} n\, \pi - \frac{15}{22}\pi) + 2\, \cos(\frac{4}{11} n\, \pi - \frac{5}{22}\pi) + 2\, \cos(\frac{6}{11} n\, \pi - \frac{3}{22}\pi) $ \\
\hphantom{$A_{11} = $} $+ 2\, \cos(\frac{8}{11} n\, \pi - \frac{3}{22}\pi) +\, 2\, \cos(\frac{10}{11} n\, \pi + \frac{5}{22}\pi).$ \\
$A_{12} = 2\, \cos(\frac{1}{6} n\, \pi - \frac{55}{72}\pi) + 2\, \cos(\frac{5}{6} n\, \pi + \frac{1}{72}\pi)$. \\
$A_{13} = 2\, \cos(\frac{2}{13} n\, \pi - \frac{11}{13}\pi) + 2\, \cos(\frac{4}{13} n\, \pi - \frac{4}{13}\pi) + 2\, \cos(\frac{6}{13} n\, \pi - \frac{1}{13}\pi) $ \\
\hphantom{$A_{13} = $} $+ 2\, \cos(\frac{8}{13} n\, \pi + \frac{1}{13}\pi) + 2\, \cos\frac{10}{13} n\, \pi + 2\, \cos(\frac{12}{13} n\, \pi + \frac{4}{13}\pi)$. \\
$A_{14} = 2\, \cos(\frac{1}{7} n\, \pi - \frac{13}{14}\pi) + 2\, \cos(\frac{3}{7} n\, \pi - \frac{3}{14}\pi) + 2\, \cos(\frac{5}{7} n\, \pi - \frac{3}{14}\pi)$. \\
$A_{15} = 2\, \cos(\frac{2}{15} n\, \pi - \frac{1}{90}\pi) + 2\, \cos(\frac{4}{15} n\, \pi - \frac{7}{18}\pi) + 2\, \cos(\frac{8}{15} n\, \pi - \frac{19}{90}\pi)$ \\
\hphantom{$A_{15} = $} $ + 2\, \cos(\frac{14}{15} n\, \pi + \frac{7}{18}\pi)$. \\
$A_{16} = 2\, \cos(\frac{1}{8} n\, \pi + \frac{29}{32}\pi) + 2\, \cos(\frac{3}{8} n\, \pi + \frac{27}{32}\pi) + 2\, \cos(\frac{5}{8} n\, \pi + \frac{5}{32}\pi) $ \\
\hphantom{$A_{16} = $} $+ 2\, \cos(\frac{7}{8} n\, \pi + \frac{3}{32}\pi)$. \\
$A_{17} = 2\cos(\frac{2}{17}n\,\pi + \frac{14}{17}\pi) + 2\cos(\frac{4}{17}n\, \pi - \frac{8}{17}\pi) + 2\cos(\frac{6}{17}n\, \pi - \frac{5}{17}\pi)$ \\
\hphantom{$A_{17} = $} $ + 2\cos\frac{8}{17}n\, \pi + 2 \cos(\frac{10}{17}n\, \pi - \frac{1}{17}\pi) + 2 \cos(\frac{12}{17}n\, \pi - \frac{5}{17} n\, \pi) $ \\
\hphantom{$A_{17} = $} $+ 2 \cos(\frac{14}{17}n\, \pi - \frac{1}{17}\pi)  + 2\cos(\frac{16}{17}n\, \pi + \frac{8}{17}\pi).$ \\
\noindent $A_{18} = 2 \cos(\frac{1}{9}n\, \pi + \frac{20}{27}\pi) + 2 \cos(\frac{5}{9}n\, \pi - \frac{2}{27}\pi) + 2 \cos(\frac{7}{9}n\, \pi + \frac{2}{27}\pi).$
\end{quotation}

What is perhaps most striking about this list is that it is entirely comprised of sums of only cosines.  And since the evaluation of Equation~\eqref{Aqn} is arguably the most complicated part of quoting Theorem~\ref{HRfull}, it is natural to ask whether this pattern indicates some simplified structure which could be exploited. 

The answer was provided by A.~Selberg, although he did not publish the result himself, as he explained in an interview~\cite{Selberg}:

\begin{quotation}\color{RedOrange}
$\cdots$  So I started to take a look at Hardy's and Ramanujan's
article on ``partitions", and I found the exact formula. But that turned out to be a
disappointment! You see, I had finished my investigation of the partition function
in the summer of 1937, and when I came to Oslo and looked in Zentralblatt, I found
the review of my first article. On the same page was the review of Rademacher's
article on the partition function. I had something that Rademacher did not have,
and that was a much simpler expression for the coefficients that appear in the
relevant series. This was undoubtedly something Ramanujan would have done if
he had been at his full power when this work was under way. $\cdots$ I made up my mind not to publish this result about the coefficients that I mentioned. I thought it was too little to write about.
\end{quotation}

This simpler expression for the coefficients was ultimately communicated to and published by A.~L.~Whiteman in \emph{A sum connected with the series for the partition function}~\cite{Whiteman}:

\begin{quotation}\color{RoyalBlue}
Because of the intricacy of the numbers $\omega_{n,k}$ the task of evaluating $A_k(n)$ for large $k$ directly from its definition in~\eqref{Aqn} is quite formidable.  To surmount this difficulty D.~H.~Lehmer~\cite{Lehmer2} made an intensive study of the $A_k(n)$.  

$\cdots$

Some years ago Atle Selberg proved (but did not publish) the result that $A_k(n)$ may be expressed alternatively in the form 
\[ A_k(n) = \left(\frac{k}{3}\right)^{1/2} \sum_{(3\ell^2+\ell)/2\equiv -n(\text{mod}\, k)} (-1)^\ell \cos\frac{6\ell+1}{6k}\ \pi,\]
where $\ell$ runs over integers in the range $0 \leq \ell < 2k$ which satisfy the summation condition.  In this striking formula $A_k(n)$ is expressed as a sum which involves only cosines and which contains considerably fewer terms than \eqref{Aqn}.  
\end{quotation}

This more efficient formulation of the $A_k(n)$ was exploited recently by F.~Johansson in~\cite{johansson2012efficient} for a nearly optimal computer implementation of $p(n)$ with running time~$O(\sqrt{n}\, \log^{4+o(1)}(n))$.


\section{What is a formula?}\label{whatisaformula}

There is a short paper by H.~Wilf titled \emph{What is an answer?}~\cite{WilfAnswer}.  
The principle observation is:

\begin{quotation}{\color{Magenta}
In many branches of pure mathematics it can be surprisingly hard to recognize when a question has, in fact, been answered.  A clearcut proof of a theorem or the discovery of a counterexample leaves no doubt in the reader's mind that a solution has been found.  But when an ``explicit solution" to a problem is given, it may happen that more work is needed to evaluate that ``solution," in a particular case, than exhaustively to examine all of the possibilities directly from the original formulation of the problem.  In such a situation, other things being equal, we may justifiably question whether the problem has in fact been solved.  

$\cdots$

For concreteness, suppose that for each integer $n>0$ there is a set $S_n$ that we want to count.  Let $f(n) = |S_n|$ (the cardinality of $S_n$), for each $n$.

Suppose further that a certain formula has been found, say
\begin{equation}\tag{W1} f(n) = \text{Formula}(n)\qquad (n=1,2,\ldots)\label{Wilf1}\end{equation}
in which Formula$(n)$ may involve various multiple summation signs extending over various sets and various complicated summands, etc.

In order to evaluate the ``answer"~\eqref{Wilf1}, let's look at the competition.  To insure my own immortality on the subject, I am now going to show you a simple formula that ``answers" all such questions at once.  If you're ready, then, here it is:
\begin{equation}\label{Wilf2}
\tag{W2} f(n) = \sum_{S_n} 1.
\end{equation}

Well, anyway, the summand is elegant, even if the range of summation is a bit untidy.

Now~\eqref{Wilf2} is unacceptable in polite society as an answer, despite its elegant appearance, because it is just a restatement of the question, and it does not give us a tool for calculating $f(n)$ that we didn't have before.  It does illustrate an important point, though: there is no counting problem for which a formula does not already exist, namely~\eqref{Wilf2}.

A first criterion for evaluating an ``answer" then, might be that ``Formula$(n)$" should be an improvement over the insightful contribution~\eqref{Wilf2}.

$\cdots$

How do we compare the complexity of evaluating Formula$(n)$ with~\eqref{Wilf2}?

$\cdots$

The functions that will be compared are, therefore:

Count$(n) = $ the complexity of the algorithm for calculating $f(n)$, whether it be given by a formula, an algorithm, et cetera, and

List$(n) = $ the complexity of producing all of the members of the set $S_n$, one at a time, by the speediest known method, and counting them.
\vskip .1in
DEFINITION 1: We will say that a solution of a counting problem is \emph{effective} if
\begin{equation*} \lim_{n\to\infty} \frac{\text{Count}(n)}{\text{List}(n)} = 0. \end{equation*}
\vskip .1in
What we are saying is that a formula or whatever is an effective solution of a problem if the effort involved in counting the members of $S_n$ by means of that formula is asymptotically small compared to the effort of constructing all of the members of $S_n$, by the best-known algorithm, and counting them.

}\end{quotation}

Equation~\eqref{HRFormula}, even before the work of Rademacher, is thus an \emph{effective} ``answer," since it implies the existence of a finite $n_0$ such that for $n \geq n_0$, we have Count$(n) = O(n^{5/2})$ (see~\cite[Section 2]{johansson2012efficient} which assumes a na\"ive implementation of Equation~\eqref{HRFormula}).
Equation~\eqref{realHR} and all its weaker corollaries, on the other hand, are not.  

This style of thought now motivates a final question: what is a formula?  Let us suppose, as in the above, that there is a set $S_n$ that we want to count, and let $f(n) = |S_n|$.  

To insure my own immortality on the subject, I am now going to demonstrate a simple formula with quantitative error rate:
\begin{equation} f(n) = 1 + O(f(n)). \label{universal} \end{equation}
Of course, insofar as the above formula is universal, it can hardly be expected to be useful.  A first criterion for evaluating a formula then, might be that it should be an improvement over the insightful contribution~\eqref{universal}.

\begin{definition}\label{defformula}
Suppose we have two formulas for counting $f(n)$, of the form
\begin{equation}\label{formula1} f(n) = g_1(n) + O(h_1(n)), \end{equation}
and 
\begin{equation}\label{formula2} f(n) = g_2(n) + O(h_2(n)), \end{equation}
for some explicit, known functions $g_1, g_2, h_1, h_2$.  
We say~\eqref{formula1} is an improvement over~\eqref{formula2} if 
\[ \lim_{n\to\infty} \frac{h_1(n)}{h_2(n)} = 0. \]
\end{definition}

In order to motivate the emphasis of the definition on the error rate rather than the estimates themselves, we now present the final chapter in this story, one of simultaneous discovery, in this case by J.~V.~Uspensky in 1920~\cite{Uspensky}.

\begin{theorem}[Uspensky, 1920~\cite{Uspensky}]\label{UspenskyTheorem}
We have\footnote{We fix a purely typographical error in the final formula where $\sqrt{\pi}$ should have been simply $\pi$. }
\begin{equation}\tag{U1}\label{U1}
p(n) = \frac{e^{\pi \sqrt{\frac{2}{3}\left(n - \frac{1}{24}\right)}}}{4\sqrt{3}\left(n-\frac{1}{24}\right)}\left(1 - \frac{\sqrt{3}}{\pi \sqrt{2n - \frac{1}{12}}}\right)  + \rho_n\, e^{\frac{\pi}{\sqrt{6}}\, \sqrt{n}},
\end{equation}
where $\rho_n$ is a function that remains bounded indefinitely growing with $n$.  
\end{theorem}

Equation~\eqref{U1} is distinct from~\eqref{realHR} only in the big-O error term; letting $h_2(n)$ denote the big-O error term in~\eqref{U1}, and $h_1(n)$ denote the big-O error term in~\eqref{realHR}, we have 
\[ \lim_{n\to\infty} \frac{h_1(n)}{h_2(n)} = \frac{1}{n} = 0. \]

This final distinction, that of a formula being defined not just in terms of its estimate of the dominant term, but also that of its error term, justifies our conclusion that~\eqref{realHR} is the real Hardy-Ramanujan formula, and all of the other Hardy-Ramanujan formulas are just imitating. 

\section{List, list, O list! [sic]}

In order to create a more perfect union of nuanced ideas into the prose of forthcoming papers, we present the following summarizing key points. 

\begin{enumerate}
\item Never mention the ``genius" of Hardy and Ramanujan in any formal work.
\item If you must use the word ``genius," never apply it to formula~\eqref{pn1}.
\item Hardy and Ramanujan's infinite series for $p(n)$ is divergent.  (Hardy and Ramanujan did not ``express any definite opinion" on the matter of its convergence or divergence, and Lehmer proved it was divergent.)
\item Rademacher's series for $p(n)$ is convergent.
\item Rademacher's convergent series for $p(n)$ has a completely effective error bound which provides explicit, numerical guarantees on the error for each $n$ and partial sum.
\item Hardy and Ramanujan themselves did not provide a means with which to compute $p(n)$ with any certainty for any specified, explicit $n$.  Rademacher used the completely effective error bounds for his series to provide completely effective error bounds for Hardy and Ramanujan's series.  
\item Lehmer used the convergence of Rademacher's series to demonstrate that for $n \geq 600$, the value of Hardy and Ramanujan's summation $\sum_{q=1}^{\frac{2}{3} \sqrt{n}} A_q \phi_q$ rounded to the nearest integer is equal to $p(n)$. 
\item Lehmer improved the error estimates on Rademacher's error term and provided an explicit parameterized theorem which, for any given value of $n$, determines how many terms you must take of Rademacher's series in order to be within a given delta of the true value of $p(n)$.
\item The first term of the Hardy-Ramanujan asymptotic formula was independently published by Uspensky several years after Hardy and Ramanujan's original paper was published.
\end{enumerate}

We next present some prose which might help alleviate future writers from inadvertently misrepresenting these various works. 

\begin{enumerate}
\item The formula~\eqref{pn1} is a vastly simplistic corollary of the original results obtained by Hardy and Ramanujan~\cite{HR}, which have much greater depth. 
\item The formula~\eqref{realHR} is the first term of an asymptotic series initially discovered by Hardy and Ramanujan~\cite{HR}; it was also independently discovered by Uspensky~\cite{Uspensky} several years later, with a slightly weaker error term.
\item Hardy and Ramanujan's asymptotic series~\cite{HR}, the first term of which is represented in~\eqref{realHR}, helped motivate subsequent work by Rademacher~\cite{rademacher1937convergent, Rademacher} in finding a convergent asymptotic series with quantitative error rate, and Lehmer~\cite{Lehmer1, Lehmer2, Lehmer} to investigate the convergence and improved error estimates of the two series.

\end{enumerate}

We end our discussion with a conjecture.  

\begin{conjecture}
There exists a proof of Proposition~\ref{propHR} via elementary methods.
\end{conjecture}

By elementary, we are explicitly excluding those methods which appeal to complex analysis, or similar analyses of comparable abstraction. 
One suggestion would be to modify Erd\H{o}s's original approach~\cite{ErdosElementary}.  Another would be to modify E.~R.~Canfield's~\cite{CanfieldElementary}.  One could also try probabilistic methods.  
The most important feature is the error term; one cannot claim to have solved the conjecture unless the error term is the same as~\eqref{realHR}.
\\

\noindent {\bf Acknowledgement.}
I am grateful to Richard Arratia for insightful discussions and feedback on earlier drafts.


\end{document}